\documentclass{amsart}

\newtheorem{theorem}{Theorem}[section]
\newtheorem*{maintheorem}{Theorem}
\newtheorem{lemma}[theorem]{Lemma}
\newtheorem{proposition}[theorem]{Proposition}
\newtheorem{corollary}[theorem]{Corollary}

\theoremstyle{definition}
\newtheorem{definition}[theorem]{Definition}

\newtheorem{remark}[theorem]{Remark}
\newtheorem*{mainremark}{Remark}
\newtheorem*{acknowledgement}{Acknowledgement}

\theoremstyle{remark}

\renewcommand{\labelenumi}{(\roman{enumi})}

\DeclareFontFamily{U}{wncy}{}
\DeclareFontShape{U}{wncy}{m}{n}{<->wncyr10}{}
\DeclareSymbolFont{mcy}{U}{wncy}{m}{n}
\DeclareMathSymbol{\Sh}{\mathord}{mcy}{"58}

\newcommand\mynote[1]{\marginpar{\ \\ \small \tt #1}}
\newcommand\bel[1]{{\mynote{#1}}\begin{equation}\label{#1}}
\newcommand\mylabel[1]{\label{#1}}

\newcommand{\ZZ}{\mathbb{Z}}
\newcommand{\QQ}{\mathbb{Q}}
\newcommand{\RR}{\mathbb{R}}
\newcommand{\CC}{\mathbb{C}}

\newcommand{\HH}{\mathbb{H}}
\newcommand{\PP}{\mathbb{P}}

\newcommand{\GG}{\mathbb{G}}

\newcommand  {\shF}     {\mathcal{F}}

\newcommand  {\shL}     {\mathcal{L}}

\newcommand  {\Aut}     {\operatorname{Aut}}

\newcommand  {\ch}      {\operatorname{ch}}

\renewcommand{\cong}    {\equiv}

\newcommand  {\End}     {\operatorname{End}}

\newcommand  {\ev}      {{\text{\rm ev}}}

\newcommand  {\Hilb}    {\operatorname{Hilb}}

\newcommand  {\id}      {\operatorname{id}}

\renewcommand  {\ker }  {\operatorname{ker}}
\newcommand  {\Km }     {\operatorname{Km}}

\newcommand  {\lra}     {\longrightarrow}

\newcommand  {\NS}      {\operatorname{NS}}

\renewcommand{\O}       {\mathcal{O}}

\newcommand  {\Pic}     {\operatorname{Pic}}

\newcommand  {\pr}      {\operatorname{pr}}

\newcommand  {\ra}      {\rightarrow}

\newcommand  {\rank}    {\operatorname{rank}}

\newcommand  {\Sym}     {\operatorname{Sym}}

\newcommand  {\Todd}    {\operatorname{Todd}}

\def\mydate{\number\day\space\ifcase\month \or January\or February\or March\or 
April\or May\or June\or July\or
August\or September\or October\or November\or December\fi \space\number\year}

\DeclareFontFamily{U}{wncy}{}
\DeclareFontShape{U}{wncy}{m}{n}{<->wncyr10}{}
\DeclareSymbolFont{mcy}{U}{wncy}{m}{n}
\DeclareMathSymbol{\Sh}{\mathord}{mcy}{"58}

\usepackage{amscd,amssymb,amsmath}
\usepackage{mathrsfs}
\usepackage[all]{xy}

\begin{document}

\title[Enriques manifolds]
      {Enriques manifolds}

\author[Keiji Oguiso]{Keiji Oguiso}
\address{Department of Mathematics, 
Osaka University, Toyonaka 560-0043 Osaka, Japan}
\curraddr{}
\email{oguiso@math.sci.osaka-u.ac.jp}

\author[Stefan Schr\"oer]{Stefan Schr\"oer}
\address{Mathematisches Institut, Heinrich-Heine-Universit\"at,
40225 D\"usseldorf, Germany}
\curraddr{}
\email{schroeer@math.uni-duesseldorf.de}

\subjclass[2000]{14J28, 14J32}

\dedicatory{Final version, February 4, 2011}

\begin{abstract}
Using the theory of hyperk\"ahler manifolds, 
we generalize the notion of Enriques surfaces to higher dimensions, explore their properties,
and construct several examples using group actions on Hilbert schemes
of points or moduli spaces of stable sheaves.
\end{abstract}

\maketitle
\tableofcontents
\renewcommand{\labelenumi}{(\roman{enumi})}

\section*{Introduction}

Naturally, Enriques surfaces play a prominent role in the Enriques classification of algebraic surfaces.
They are, by definition,  minimal surfaces of Kodaira dimension $\kappa=0$ with $b_2=10$. An equivalent condition
is that they are not simply connected and have a K3 surface as universal covering.
There are numerous deep results concerning Enriques surfaces, for example about their geometry,  automorphism groups,
or periods. In light of this richness it is natural to ask whether there is a natural generalization of
 Enriques surfaces to higher dimensions.
The goal of this paper is to introduce the notion of \emph{Enriques manifolds}, explore their basic properties,
and construct several interesting examples.

Recall that a \emph{hyperk\"ahler manifold} is a smooth compact simply-connected K\"ahler manifold $X$ with the property
that $H^0(X,\Omega^2_X)$ is generated by a symplectic form. Beauville \cite{Beauville 1983} showed that such manifolds are, together
with complex tori and Calabi--Yau manifolds, the basic building blocks for K\"ahler manifolds with $c_1=0$.
There is a profound theory for hyperk\"ahler manifolds (see Huybrechts \cite{Huybrechts 1999}), which largely runs parallel to 
the theory of K3 surfaces. Indeed, one should view hyperk\"ahler manifolds as the correct generalization of  
K3 surfaces to higher dimensions. Therefore, we define an \emph{Enriques manifold} as a  connected complex space $Y$ that is not
simply connected and whose universal covering $X=\widetilde{Y}$ is a hyperk\"ahler manifold.

It turns out that the fundamental group $\pi_1(Y)$ is a finite cyclic group. We call its order $d\geq 2$  the \emph{index}
of the Enriques manifold $Y$. This number $d$ is  a divisor of $n+1$, where $\dim(Y)=2n$, and moreover
meets the condition $\varphi(d)<b_2(X)$, where $\varphi$ is Euler's phi function. A natural question arises: Which
integers $d$ appear as indices for Enriques manifolds? In other words, which cyclic groups can act freely on
some hyperk\"ahler manifold? 

In some sense, there are not too many known examples of hyperk\"ahler manifolds. Beauville \cite{Beauville 1983} constructed two
infinite series, namely the Hilbert scheme of points for K3 surfaces,
and his \emph{generalized Kummer variety} $\Km^n(A)$, defined as a   Bogomolov factor in $\Hilb^{n+1}(A)$ for abelian surfaces $A$. 
Furthermore, there are two sporadic examples of O'Grady (see \cite{O'Grady 1999}, \cite{O'Grady 2003}). The first idea to construct Enriques manifolds is
to look at Hilbert schemes for Enriques surfaces or bielliptic surfaces, but this does not work out. Rather, 
it leads to an interesting new construction of \emph{Calabi--Yau manifolds}:

\begin{maintheorem}
Let $S$ be an Enriques surface or a bielliptic surface, and $n\geq 2$. Then $\Hilb^n(S)$ has a finite \'etale covering
that is a Calabi--Yau manifold or is the product of  a Calabi--Yau manifold with an elliptic curve, respectively.
\end{maintheorem}

However, if one starts with an Enriques surface $S'$, say with universal covering $S$, and an \emph{odd} number $n\geq 1$, then the induced action
of $G=\pi_1(S')$ on $X=\Hilb^n(S)$ is free, and the corresponding quotient is an Enriques manifold $Y$ 
of dimension $\dim(Y)=2n$ and index $d=2$. There is   a variant with generalized Kummer varieties,
and the preceding construction can be extended from Hilbert schemes of points to moduli spaces of sheaves:

\begin{maintheorem}
Suppose $S'$ is an Enriques surface whose corresponding K3 surface has Picard number $\rho(S)=10$.
Let $v=(r,l,\chi-r)\in H^\ev(S,\ZZ)$ be a primitive Mukai vector with $v^2\geq 0$
and $\chi\in\ZZ$ odd. Then for very general polarizations $H\in\NS(S)_\RR$, 
the moduli space $X=M_H(v)$ is a hyperk\"ahler manifold endowed with a free action of $G=\pi_1(S')$, and 
$Y=X/G$ is an Enriques manifold of dimension $v^2+2$ and index $d=2$.
\end{maintheorem}

Recall that a \emph{bielliptic surface} $S$ has, by definition, a finite \'etale covering $A$ that is an abelian surface.
To construct examples of Enriques of higher index, we use the classification of bielliptic surface due to
Bagnera and de Franchis and study the induced action of $G=\pi_1(S)$ on $\Hilb^n(A)$. This yields:

\begin{maintheorem}
There are Enriques manifolds with index $d=2,3,4$.
\end{maintheorem}

The paper is organized as follows:
In the first section we recall several results about the Bogomolov decomposition of manifolds
with trivial first Chern class and the theory of hyperk\"ahler manifolds.
In the second section, we introduce the notion of Enriques manifolds and collect their basic
properties. In the third section, we examine Hilbert schemes of points for Enriques surfaces
and bielliptic surfaces. The first examples of Enriques manifolds appear in Section 4
as quotients of Hilbert schemes of points for the K3 covering of an Enriques surface.
We extend this construction to moduli spaces of stable sheaves in Section 5.
In Section 6 we use the classification of bielliptic surfaces to construct Enriques manifolds whose
universal covering are Beauville's generalized Kummer varieties.

\begin{mainremark}
Our paper overlaps  with the paper  \emph{Higher dimensional Enriques varieties and automorphisms of generalized Kummer varieties},
arXiv:1001.4728v1 by Samuel Boissi\`ere, Marc Nieper-Wisskirchen, Alessandra Sarti, which appeared on the arXiv after our paper was completed.
\end{mainremark}

\begin{acknowledgement}
The authors started this research at the Centro di Ricerca Matematica Ennio De Giorgi during the workshop
on \emph{Groups in Geometry} in September 2008. We thank the Centro for providing a stimulating research environment.
The second author visited the  first author in September 2009 at the Osaka University and thanks the Department of Mathematics for
its hospitality. The second author was supported by a DFG grant within the Forschergruppe FOR 790
\emph{Classification of Algebraic Surfaces and Compact Complex Manifolds}. Finally, we thank the referee, Alessandra Sarti and
Samuel Boissi\'ere for   comments.
\end{acknowledgement}

\section{Bogomolov decomposition and hyperk\"ahler manifolds}
\mylabel{bogomolov decomposition}

Throughout this paper, we shall work over the complex numbers.
Given a complex manifold  $Y$, we regard its first Chern class
$c_1(Y)$ as an element in the rational vector space $H^2(Y,\QQ)$.
In this section we recall some results on  compact K\"ahler manifolds $Y$ with $c_1(Y)=0$,
which are due to Beauville,  Bogomolov, Fujiki, Huybrechts, Mukai, O'Grady, Yoshioka, and others.

The fundamental result is that such manifolds $Y$ admit a finite \'etale covering 
$X\ra Y$ of the form  $X=\prod_{i=1}^r X_i$
where the factors $X_i$ are complex tori, Calabi--Yau manifolds, or hyperk\"ahler manifolds
(\cite{Beauville 1983} and \cite{Bogomolov 1974b}).
Such a factorization on a finite \'etale cover is called a \emph{Bogomolov decomposition}.
The fundamental group of $Y$ is an extension of a finite group by a free abelian group.
Obviously, $\pi_1(Y)$  is finite if and only if no Bogomolov factor is a complex torus.
In this case, a Bogomolov factorization exists only on the universal covering $X=\widetilde{Y}$.

Throughout the paper, the term  \emph{Calabi--Yau manifold} denotes a compact connected K\"ahler
manifold $X$ of dimension $\geq 3$ that 
is simply connected, has $\omega_X=\O_X$, and $h^{p,0}(X)=0$ for $0<p<\dim(X)$. 
With this definition, Calabi--Yau manifolds are automatically projective, by Kodaira's embedding Theorem.
There is no common agreement about   the term ``Calabi--Yau manifold'',
and some authors use it to denote manifolds with $c_1=0$. Also, it is sometimes useful to
replace the assumption that $X$ is K\"ahler by the weaker assumption that $X$ is is bimeromorphic
to some K\"ahler manifold, which then is equivalent to being Moishezon.

Recall that a \emph{hyperk\"ahler manifold} is a compact connected K\"ahler manifold $X$
that is simply-connected with $H^0(X,\Omega^2_X)=\CC\sigma$, where $\sigma$ is a symplectic form.
In other words, it induces nondegenerate alternating pairing on all tangent spaces $\Theta_X(x)$, $x\in X$. 
Let us recall some facts on such manifolds.
The existence of a  symplectic form $\sigma\in H^0(X,\Omega^2_X)$ ensures that
$\dim(X)=2n$ is even and that the dualizing sheaf $\omega_X=\O_X$ is trivial.
Moreover, one knows that the algebra of holomorphic forms $\bigoplus_p H^0(X,\Omega^p_X)$
is generated by the symplectic form, such that
$$
h^{p,0}(X)=
\begin{cases}
1 & \text{if $p$ is even,}\\
0 & \text{if $p$ is odd,}
\end{cases}
$$
thus $\chi(\O_X)=n+1$.
There is an elaborate general theory about hyperk\"ahler manifolds 
parallel to the theory of K3 surfaces, see Huybrechts  \cite{Huybrechts 1999}.

Up to deformation equivalence, the only known examples of hyperk\"ahler manifolds are
the Hilbert scheme of points $\Hilb^n(S)$ for K3 surfaces $S$ and Beauville's generalized Kummer varieties $\Km^n(A)$ for abelian
surfaces $A$, both introduced by Beauville \cite{Beauville 1983}, 
and two sporadic examples $M_6$, $M_{10}$ constructed by O'Grady as desingularizations of certain moduli spaces
of   sheaves on K3 or abelian surfaces (\cite{O'Grady 1999} and \cite{O'Grady 2003}). Here are some numerical invariants for these hyperk\"ahler manifolds:
$$
\renewcommand{\arraystretch}{1.5}
\begin{array}[t]{c|c|c|c|c}
X          &\Hilb^n(S)  & \Km^n(A)     & M_6  & M_{10}\\
\hline
\dim(X)    & 2n         & 2n           &    6 & 10\\
\hline
\chi(\O_X) & n+1        & n+1          &    4 & 6\\
\hline
b_2(X)     & 23         & 7            &    8 & 24 \\
\end{array}
$$

Next, let us recall some  facts about symmetric products and  Hilbert schemes of points.
Given an arbitrary compact complex space $V$, we denote
by $\Sym^n(V)$ the symmetric product. Let $\pi:V^n\ra\Sym^n(V)$ be the quotient map.
If $V$ is normal, so is $\Sym^n(V)$, and the norm map $\pi_*(\O_{V^n}^\times)\ra\O_{\Sym^n(V)}^\times$
induces a homomorphism
$$
\Pic(V)\lra\Pic(\Sym^n(V)),\quad \shL\longmapsto\shL^{(n)}.
$$
Using the interpretation $\shL^{(n)}=\pi_*(\otimes_{i=1}^n\pr_i^*(\shL))^{S_n}$, we easily infer
$H^0(V,\shL)=H^0(\Sym^n(V),\shL^{(n)})$, such that the preceding
homomorphism is injective. If $V$ is Gorenstein, so is  $\Sym^n(V)$, and we have
$\omega_{\Sym^n(V)}=(\omega_V)^{(n)}$.

According to Grothendieck's observation (\cite{SGA 1}, Expose IX, Remark 5.8), we have an identification
$$
\pi_1(\Sym^n(V)) = H_1(V,\ZZ),\quad n\geq 2.
$$
If $V$ is normal with only quotient singularities, then $\Sym^n(V)$ is normal with  only quotient singularities.
Under this assumption, the canonical map
$$
\pi_1(Z)\lra\pi_1(\Sym^n(V))=H_1(V,\ZZ)
$$
is bijective for any resolution of singularities $Z\ra\Sym^n(V)$, by \cite{Kollar 1993}, Theorem 7.8. 
Now consider the Hilbert scheme, or rather Douady space, of points $\Hilb^n(V)$, and let
$$
\gamma:\Hilb^n(V)\lra\Sym^n(V),\quad A\longmapsto\sum_{x\in V}\operatorname{length}(\O_{A,x})x
$$
be the \emph{Hilbert--Chow morphism}, which sends a subscheme   to the corresponding zero-cycle 
(see \cite{Iversen 1970} for more details).
In general, the Hilbert scheme of points is much more complicated than  the symmetric product.
However, if $V$ is a smooth surface, then  the Hilbert--Chow morphism  
is a crepant resolution of singularities. 
We refer to  Beauville's paper \cite{Beauville 1983} or the monograph of Brion and Kumar \cite{Brion; Kumar 2005}, Chapter  7 
for   detailed discussions.

Now let us recall Beauville's generalized Kummer surface.
Let  $A$ be an abelian surface, and consider the composite map
$$
\Hilb^{n+1}(A)\lra\Sym^{n+1}(A)\lra A,
$$
where the first arrow is the Hilbert--Chow morphism, and the second arrow is the addition map.
This actually is the Albanese map. The fiber over the origin $\Km^n(A)\subset\Hilb^{n+1}(A)$ is called the 
\emph{generalized Kummer variety}, and is a hyperk\"ahler manifold $X=\Km^n(A)$ of dimension $2n$.
In other words, $\Km^n(A)$ is defined as a Bogomolov factor for $\Hilb^{n+1}(A)$.

Moduli spaces of coherent sheaves provide further examples.
Let $S$ be a K3 surface. $v\in H^\ev(S,\ZZ)$   a Mukai vector, and $H\in\NS(S)_\RR$
a polarization. Mukai \cite{Mukai 1984} showed that the moduli space $M_H(v)$ of $H$-stable sheaves $\shF$ on $S$ with
Mukai vector $v(\shF)=v$ is smooth of dimension $v^2+2$, where   $v^2=(v,v)$ comes from the 
Mukai pairing (for details, see   \cite{Mukai 1987} or \cite{Huybrechts; Lehn 1997}, Chapter 6).
It turns out that for $H$ generic and $v$ primitive, $M_H(v)$ is actually a hyperk\"ahler 
manifold (see \cite{O'Grady 1997} and \cite{Yoshioka 2001}), which is deformation equivalent to $\Hilb^n(S)$, $n=(v^2+2)/2$.
Using moduli spaces of stable sheaves on abelian surfaces and the Fourier--Mukai transform,
Yoshioka \cite{Yoshioka 2001} constructed hyperk\"ahler manifolds $K_H(v)$ that are deformation equivalent to $\Km^n(A)$.

There are more examples of hyperk\"ahler $4$-folds,   all of them
deformation equivalent to   $\Hilb^2(S)$:
Beauville and Donagi \cite{Beauville; Donagi 1985} showed that the variety of lines on a smooth
cubic hyperplanes in $\PP^5$ is a hyperk\"ahler $4$-fold.
Iliev and Ranestad \cite{Iliev; Ranestad 2001} proved that the variety of sums of powers for a general
cubic hyperplane as above is an another such examples.
O'Grady \cite{O'Grady 2006} constructed hyperk\"ahler $4$-folds as double covers of
certain sextic hyperplane in $\PP^5$.
Debarre and Voisin \cite{Debarre; Voisin 2009} showed that for $V=\CC^{\oplus 10}$ and $\sigma\in\Lambda^3(V^\vee)$ general,
the scheme of $6$-dimensional subvector spaces of $V$ on which $\sigma $ vanishes is
a hyperk\"ahler $4$-fold.

\section{Notion of Enriques manifolds}
\mylabel{notion enriques}

In the classification of surfaces,  \emph{Enriques surfaces} are defined
as minimal surfaces $S$ with Kodaira dimension $\kappa=0$ and second Betti number $b_2=10$.
A different but equivalent definition is that $S$ is not simply connected, 
and its   universal cover  is a K3 surface.
Viewing hyperk\"ahler manifolds as the correct generalization of K3 surfaces,
we propose the following generalization of Enriques surfaces to higher dimensions. 

\begin{definition}
\mylabel{enriques manifolds}
An \emph{Enriques manifold} is a connected complex manifold $Y$
that is not simply connected and whose universal cover $X$ is 
a hyperk\"ahler manifold.
\end{definition}

Obviously, Enriques manifolds $Y$ are compact, 
of even dimension $\dim(Y)=2n$,
and with finite fundamental group. Averaging a K\"ahler metric on $X$ over its $G$-translates,
one   sees that $Y$ is a K\"ahler manifold. The $2$-dimensional Enriques manifolds
are precisely the Enriques surfaces, whose fundamental group is cyclic of order $d=2$.
In higher dimensions, this order is a basic numerical invariant:

\begin{definition}
\mylabel{index}
The \emph{index} of an Enriques manifold $Y$ is the order $d\geq 2$
of its fundamental group $\pi_1(Y)$.
\end{definition}

Let $Y$ be an Enriques manifold of dimension $\dim(Y)=2n$, and $X\ra Y$ be its universal covering.
Fix a base point $y\in Y$. The natural action of 
$\pi_1(Y,y)$ on $X$ induces a representation on the 1-dimensional
vector space $H^0(X,\Omega_X^2)$, which corresponds via the trace to 
a homomorphism $\rho:\pi_1(Y,y)\ra\CC^\times$.

\begin{lemma}
\mylabel{trace}
The homomorphism of groups $\rho:\pi_1(Y,y)\ra\CC^\times$ is injective.
\end{lemma}

\proof
Let $G\subset\ker(\rho)$ be a cyclic subgroup, say of order $m=|G|$,
and consider the complex manifold $Z=X/G$.
Let $f:X\ra Z$ be the canonical projection. 
Then the $\O_Z$-algebra $f_*(\O_X)$ takes the form
$\O_Z\oplus\shL\oplus\shL^{\otimes 2}\oplus\ldots\oplus\shL^{\otimes(m-1)}$
for some $\shL\in\Pic(Z)$ of order $m$, where the multiplication   is
given by some trivialization $\shL^{\otimes m}\ra\O_Z$.
We have $\chi(\shL^{\otimes i})=\chi(\O_Z)$ because $\shL$ is numerically trivial,
whence $\chi(\O_X)=|G|\chi(\O_Z)$.

On the other hand, we have
$H^0(Z,\shF)=H^0(X,f^*(\shF))^G$ for every coherent sheaf $\shF$ on $Z$.
The canonical map $f^*(\Omega^1_Z)\ra\Omega^1_X$ is bijective, because $f:X\ra Z$ is \'etale,
and consequently $H^0(Z,\Omega_Z^p)=H^0(X,\Omega_X^p)^G$ for all $p\geq 0$.
The group $H^0(X,\Omega^p_X)$ vanishes for $p$ odd, and is 
generated by the $p$-form $\sigma\wedge\ldots\wedge\sigma$ for $p$ even.
Using that $\sigma$ is  $G$-invariant,
we conclude that the canonical maps $H^0(Z,\Omega^p_Z)\ra H^0(X,\Omega^p_X)$ are
bijective. Hodge symmetry yields
$$
\chi(\O_X)=\sum_p(-1)^p h^{p,0}(X)=\sum_p(-1)^ph^{p,0}(Z)=\chi(\O_Z),
$$
and this number equals $n+1\neq 0$. Whence $|G|=1$, and the representation $\rho $ is faithful.
\qed

\medskip

\begin{proposition}
\mylabel{index divides}
Let $Y$ be an Enriques manifold of dimension $\dim(Y)=2n$.
Then $\pi_1(Y)$ is a  cyclic group whose order $d\geq 2$ is
a divisor of $n+1$.
\end{proposition}

\proof
Being a finite subgroup of the multiplicative group of a field, 
the fundamental group must be cyclic. The universal covering $X$ has $\chi(\O_X)=n+1$, 
and we have $\chi(\O_X)=d\chi(\O_Y)$.
\qed

\begin{remark}
The representation on $H^0(X,\Omega^2_X)$ yields a canonical identification $\pi_1(Y)=\mu_d(\CC)$
with the group of $d$-th roots of unity. Consequently  we have a canonical generator $e^{2\pi\sqrt{-1}/d}$
of the fundamental group.
\end{remark}

\medskip
In the following, $Y$ denotes an Enriques manifold of
dimension $\dim(Y)=2n$ and index $d\geq 2$, with universal cover $X$.
Index and dimension control the Hodge numbers
$h^{p,q}=\dim H^q(Y,\Omega_Y^p)$ for $q=0$ and $p=0$:

\begin{proposition}
\mylabel{hodge numbers}
We have
$$
h^{0,p}(Y)=h^{p,0}(Y)=\begin{cases}
1 & \text{if $2d\mid p$ and $p\leq 2n$,}\\
0 & \text{else.}
\end{cases}
$$
In particular, $\chi(\O_Y)=(\dim(Y)+2)/2d$.
\end{proposition}

\proof
Using  $H^0(Y,\Omega_Y^p)=H^0(X,\Omega_X^p)^G$, we have $h^{p,0}(X)=0$ for $p$ odd or $p>2n$.
Consider the case $p\leq 2n$ even. Then the canonical map 
\begin{equation}
\label{tensor power}
\bigotimes_{i=1}^{p/2} H^0(X,\Omega_X^2)\lra H^0(X,\Omega_X^p),\quad
\sigma_1\otimes\ldots\otimes\sigma_{p/2}\longmapsto \sigma_1\wedge\ldots\wedge\sigma_{p/2}
\end{equation}
is bijective, and $H^0(X,\Omega_X^2)$ is $1$-dimensional. The statement now follows from Lemma \ref{trace},
together with   Hodge symmetry.
\qed

\begin{corollary}
Every Enriques manifold $Y$ is projective, and the same holds for its universal cover $X$.
\end{corollary}

\proof
The  inclusion $H^{1,1}(Y)_\RR\subset H^2(Y,\RR)$ is bijective, since $h^{2,0}=h^{0,2}=0$.
Using that the K\"ahler cone inside $H^{1,1}(Y)_\RR$ is nonempty and open, we conclude that there
must be an integral K\"ahler class on $Y$. According to Kodaira's Embedding Theorem
(compare \cite{Griffiths; Harris 1978}, p.\ 191), the K\"ahler manifold $Y$ is projective.
Pulling back this integral K\"ahler class, we obtain an integral K\"ahler class on $X$,
so $X$ is projective as well.
\qed

\begin{proposition}
The group $\Pic(Y)$ is finitely generated. Its
torsion subgroup  is a cyclic group of order $d$,
which is generated by the canonical class $\omega_Y\in\Pic(Y)$.
\end{proposition}

\proof
The map $f^{-1}(\O_Y^\times)\ra\O_X^\times$
is bijective because $f$ is \'etale. Whence we have a spectral sequence
$H^r(G,H^s(X,\GG_m))\Rightarrow H^{r+s}(Y,\GG_m)$, which  yields an exact sequence
$$
0\lra H^1(G,\CC^\times)\lra \Pic(Y)\lra\Pic(X)^G,
$$
where $G$ acts trivially on $\CC^\times$. In turn, the group cohomology
$H^1(G,\CC^\times)$ is the kernel of the map $\CC^\times\ra\CC^\times$, $z\mapsto z^d$, which
is cyclic of order $d$.
We have $\Pic^0(X)=0$ because its tangent space $H^1(X,\O_X)$ vanishes.
It follows that $\Pic(X)$ is finitely generated. Its torsion part vanishes because
$X$ is simply connected. Hence $\Pic(Y)$ is finitely generated, and its torsion part
is cyclic of order $d$.

The canonical class $\omega_Y\in\Pic(Y)$ has finite order because $f^*(\omega_Y)=\omega_X$ is trivial.
Let $r\mid d$ be its order. Then the induced $G$-action on $\omega_X^{\otimes r}=f^*(\omega_Y^{\otimes r})$
is trivial. On the other hand,  the representation on $H^0(X,\omega_X^{\otimes r})$ has trace $\tau^{rn}:G\ra\CC^\times$, which follows from
(\ref{tensor power}).
Here  $\tau$ denotes the trace of the representation on $H^0(X,\Omega^2_X)$.
Therefore $d\mid rn$. Proposition \ref{index divides} implies that $d\mid r$, whence $d=r$.
\qed

\medskip
There is a strong relation between the index of $Y$ and
the second Betti number of $X$. Let
$\varphi(d)$ be the order of the multiplicative group $(\ZZ/d\ZZ)^\times$.
Recall that $\varphi(d)=\prod p_i^{\nu_i-1}(p_i-1)$  if $d=\prod_i p_i^{\nu_i}$ is the prime
factorization.

\begin{proposition}
We have $\varphi(d)<b_2(X)$.
\end{proposition}

\proof
This result is due to Nikulin for K3-surfaces  (\cite{Nikulin 1980}, Theorem 3.1); his argument works in our situation
as follows:
First note that the restriction of the Beauville--Bogomolov form $q_X$ to $H^{1,1}(X)_{\RR}\subset H^2(X,\RR)$
has index $(1,b_2(X)-3)$, and that $q_X(\alpha)>0$ for all K\"ahler classes $\alpha\in H^2(X,\RR)$, as explained
in \cite{Huybrechts 1999}, Section 1.9. Using that $X$ is projective, we infer that $q_X$ is nondegenerate on $\NS(X)\subset H^2(X,\RR)$.
Let $T\subset H^2(X,\ZZ)$ be the orthogonal complement of $\NS(X)\subset H^2(X,\ZZ)$ with
respect to the Beauville--Bogomolov bilinear form $q_X$, which clearly is a lattice of rank $<b_2(X)$.
Choose a generator $g_0\in G$. By the Cayley--Hamilton Theorem, it suffices to check that the minimal polynomial for the endomorphism
$g_0\in\End(T_\QQ)$ is the $d$-th cyclotomic polynomial $\Phi_d$, which has degree $\varphi(d)$.
Suppose this is not the case. Then there is a nontrivial $g\in G$ admitting an eigenvector $x\in T_\QQ$
with eigenvalue $1$.
Choose a generator $\sigma\in H^{2,0}(X)$ and write $g^*(\sigma)=\xi\sigma$ for some nontrivial $\xi\in\CC^\times$.
The $G$-invariance of the Beauville--Bogomolov form yields
$$
q_X(x,\sigma)=q_X(g^*(x),g^*(\sigma))=q_X(x,\xi\sigma)=\xi q_X(x,\sigma).
$$
Since $\NS(X)\subset H^2(X,\ZZ)$  
is the orthogonal complement of $\sigma$, we have $q_X(x,\sigma)\neq 0$,
whence $\xi=1$, contradiction.
\qed

\medskip
Let us tabulate the possible indices for the Betti numbers for
Beauville's families of hyperk\"ahler manifolds:
$$
\renewcommand{\arraystretch}{1.5}
\begin{array}[t]{c|c|l}
b_2 & X          & \text{possible indices $d$}\\
\hline
  7 & \Km^n(A)   & 2-10,12,14,18,24\\
\hline
 23 & \Hilb^n(S) & 2-28, 30,32,33,34,36,38,40,42,44,46,50,54,66\\
\end{array}
$$
Of course, Proposition \ref{index divides} gives stronger restrictions if $n$  has few divisors.
For example, with $n=2$ or $n=3$  we only have the possibilities
$d=3$ or $d=2,4$, respectively.
With O'Grady's 6-dimensional example $M_6$  only $d=2,4$ are possible, and
with his  10-dimensional example $M_{10}$ only $d=2,3,6$ may occur.

\medskip
We finally touch upon the subject of birationally equivalent Enriques manifolds:

\begin{proposition}
Let $Y$ and $Y'$ be two Enriques manifolds that are birationally equivalent.
Then they have the same index.
\end{proposition}

\proof
The fundamental group is a birational invariant for smooth complex manifolds
(compare \cite{Griffiths; Harris 1978}, Section 4.2), whence $d(Y)=d(Y')$.
\qed

\section{Calabi-Yau manifolds via Hilbert schemes of points}
\mylabel{hilbert schemes}

The prime goal of this paper is to construct examples of Enriques manifolds.
The first idea that comes to mind is to look at Hilbert schemes
of Enriques surfaces. This, however, does not lead to Enriques manifolds:

\begin{theorem}
\mylabel{calabi-yau}
Let $S$ be an Enriques surface and  $Y=\Hilb^n(S)$ for some $n\geq 2$.
Then $\pi_1(Y)$ is cyclic of order two, and the
universal covering $X$ of $Y$ is a Calabi--Yau manifold. 
\end{theorem}

\proof
We have $h^p(\O_S)=0$ for all $p>0$. Let $\pr:S^n\ra S^{n-1}$ be a projection.
Then $R^p\pr_*(\O_{S^n})=0$ for all $p>0$.
Induction on $n$ and the Leray--Serre spectral sequence yields $h^p(\O_{S^{n}})=0$ for all $p>0$.
Now consider the canonical projection $g:S^n\ra \Sym^n(S)$. 
The resulting  inclusion $\O_{\Sym^n(S)}\subset g_*(\O_{S^n})$ is the inclusion of invariants
with respect to the permutation action of the symmetric group, whence is  a direct summand.
It follows that $h^p(\O_{\Sym^n(S)})=0$ for all $p>0$. In turn, we have $h^p(\O_Y)=0$ for all $p>0$ because the singularities
on the symmetric product are rational. In particular, $\chi(\O_Y)=1$.

As discussed in Section \ref{bogomolov decomposition}, the dualizing sheaf $\omega_Y$ has order two,
and the fundamental group $\pi_1(Y)$ is cyclic of order two.
Taking $\dim(Y)=2n>2$ into account,
we deduce the assertion from the following Lemma.
\qed

\begin{lemma}
\mylabel{k3 or cy}
Let $Y$ be a compact K\"ahler manifold with $c_1(Y)=0$. Suppose $\chi(\O_Y)=1$
and $\pi_1(Y)$ is cyclic of order two. Then the universal covering of $Y$
is either a K3 surface or a Calabi--Yau manifold of even dimension.
\end{lemma}

\proof
Let $X\ra Y$ be the universal covering. The manifold $X$ is compact because $\pi_1(Y)$ is finite.
Consider the Bogomolov decomposition $X=\prod_{i=1}^r X_i$, where the factors are  
Calabi--Yau manifolds or hyperk\"ahler manifolds. Complex tori do not appear,
because $Y$ has finite fundamental group. We have
$$
2=|\pi_1(Y)|\chi(\O_Y)=\chi(\O_X)=\prod_{i=1}^r\chi(\O_{X_i}).
$$
In particular, $\chi(\O_{X_i})\neq 0$, so no Bogomolov factor is an odd-dimensional Calabi-Yau manifold.
Consequently $\chi(\O_{X_i})\geq 2$, which gives the estimate $2\geq 2^r$, whence $r=1$.
Thus $X=X_1$ is either an even-dimensional Calabi-Yau manifold, or 
hyperk\"ahler. In the latter case, we have $\dim(X)=2m$ and $\chi(\O_X)=m+1$,  so $m=1$ and $X$ is a K3 surface.
\qed

\medskip
What are the Bogomolov factors for the Hilbert scheme of points for bielliptic surfaces?
In this situation,  Calabi--Yau manifolds of odd dimension show up. In contrast to the even-dimensional case,
the numbers $\chi(\O_V)$ are then not so helpful. Rather than working with Euler characteristics,
we shall work with graded algebras. Suppose $V$ is a compact   K\"ahler manifold.
Then we have the cohomology algebra
$$
H^\bullet(V,\CC)=\bigoplus_iH^i(V,\CC)=\bigoplus_{p,q} H^q(X,\Omega^p),
$$
which contains the algebra of holomorphic forms
$$
R^\bullet(V,\CC)=\bigoplus_p H^0(V,\Omega^p)\subset H^\bullet(V,\CC)
$$
as a subalgebra. Note that these algebras are graded-commutative, and
that $\chi(\O_V)=\sum_p(-1)^ph^{p,0}(V)$. We shall use the following
well-known properties of the algebra of holomorphic forms:

\begin{lemma}
\mylabel{pullback injective}
Let $f:V'\ra V$ be a proper dominant morphism of compact K\"ahler manifolds.
Then the pullback map $R^\bullet(V)\ra R^\bullet(V')$ is injective.
It is even bijective provided $f$ is bimeromorphic.
\end{lemma}

%

\medskip
The following observation will allow us to obtain some information on the Bogomolov factors
in the Hilbert scheme of points for bielliptic surfaces:

\begin{proposition}
\mylabel{odd cy}
Let $X$ and $Y$ be compact connected K\"ahler manifolds of dimension $2n+1$
with $c_1=0$. Suppose that there is an finite \'etale covering $f:\tilde{X}\ra X$
so that $\tilde{X}$ is the product of an elliptic curve and a hyperk\"ahler manifold.
Assume there is a rational dominant map $r:X\dasharrow Y$ and that $\pi_1(Y)$ is finite.
Then the universal covering $\tilde{Y}$ is   a Calabi--Yau manifold.
\end{proposition}

\proof
We first reduce to the case $\pi_1(Y)=0$.
Consider the following commutative diagram:
$$
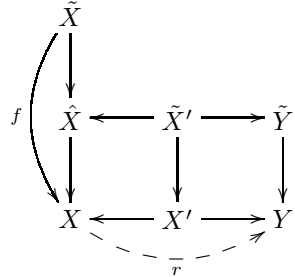
\begin{xy}
\xymatrix{
\tilde{X}\ar[d]\ar@/_1.5em/[dd]_f\\
\hat{X}\ar[d]  & \tilde{X}'\ar[l]\ar[r]\ar[d]  &  \tilde{Y}\ar[d]\\
X \ar@/_1.5em/@{-->}[rr]_r           & X'\ar[l]\ar[r]                & Y\\
}
\end{xy}
$$
Here $\tilde{Y}\ra Y$ is the universal covering, $X'$ is a smooth compact K\"ahler manifold,
$X'\ra X$ is a bimeromorphic proper  morphism, and $X'\ra Y$ is a dominant proper morphism.
The square to the right is cartesian, such that $\tilde{X}'\ra X'$ is   finite \'etale.
By bimeromorphic invariance of the fundamental group, the induced map $\pi_1(X')\ra\pi_1(X)$ is bijective,
hence there is a finite \'etale covering $\hat{X}\ra X$ making the square to the left cartesian.
We may also assume that our given finite \'etale covering $\tilde{X}\ra X$ factors over $\hat{X}$,
by passing to a larger covering.
Replacing $Y,X$ by $\widetilde{Y},\widetilde{X}$, we may assume that $Y$ is simply connected
and that $X=M\times E$ is the product of a hyperk\"ahler manifold $M$
and an elliptic curve $E$.

The idea now is to use the algebra of holomorphic forms
$$
R^\bullet(V)=\bigoplus_{p\ge 0}H^0(V,\Omega^p_V)
$$
for various schemes $V$ occurring in our situation. 
Let $Y=Y_1\times\ldots\times Y_d$ be the Bogomolov decomposition.
There are no abelian factors because $\pi_1(Y)=0$.
Consider the chain of injective pull back maps 
$$
R^\bullet(X)\lra R^\bullet(X')\longleftarrow R^\bullet(Y)
$$
induced from the preceding diagram. The map on the left is bijective since the proper morphism
$X'\ra X$ is bimeromorphic, by Lemma \ref{pullback injective}. 
Therefore, we may regard $R^\bullet(Y)\subset R^\bullet(X)$
as a  subalgebra.  Moreover, $R^\bullet(X)=R^\bullet(M)\otimes R^\bullet(E)$.
Note that the product of elements of odd degree is zero, because $\dim(E)=1$.

Seeking a contradiction, we now suppose that one Bogomolov factor $Y_i$ is hyperk\"ahler, say
of dimension $2n_i\leq 2n-2$. We infer
that there is a nonzero element  $\sigma\in R^\bullet(Y)$ of degree two with $\sigma^{n_i+1}=0$.
On the other hand, $R^\bullet(M)$ is necessarily generated by $\sigma$, and therefore $\sigma^{n_i+1}\neq 0$,
contradiction. Whence all Bogomolov factors $Y_i$ are  Calabi--Yau manifolds, say of
dimension $\dim(Y_i)=n_i$, such that $2n+1=n_1+\ldots+n_d$.
Let $\sigma_i\in R^\bullet(Y_i)$ be a nonzero element of degree $n_i$. Note that $\sigma_i^2=0$,
and $\sigma_i\sigma_j\neq 0$ for $i\neq j$.

Now suppose there are two odd-dimensional Bogomolov factors, say $Y_1,Y_2$.
Then $\sigma_1\sigma_2\neq 0$ inside $R^\bullet(Y)$. But since the degrees of the  pullbacks
$r^*(\sigma_i)\in R^\bullet(X)$ are odd, we have $r^*(\sigma_1)r^*(\sigma_2)=0$, contradiction.
Whence there is at most one odd-dimensional Bogomolov factor. Since $\dim(Y)$ is odd,
there must be precisely one such factor.

It remains to show that there are no even-dimensional Bogomolov factors. 
Suppose there is such a factor, say $Y_1$. Set $m=n_1/2$.
Clearly, $m\leq n-1$. 
In order to work with functions fields, let us now assume that $X$ is a algebraic
(In the general case, one has to work with general points rather then function fields).
Choose  transcendence basis $t_1,\ldots ,t_{2m}\in K(Y_1)$ in the function field of $Y_1$.
Using the canonical projection, we may regard  $K(Y_1)\subset K(Y)$ as a subfield
of the function field of $Y$. Then the differentials $dt_1,\ldots,dt_{2m}\in\Omega^1_{K(Y)/\CC}$ are
linearly independent. Choose a transcendence basis $s_1,\ldots,s_{2n}\in K(M)$
and $s_{2n+1}\in K(E)$, such that $ds_1,\ldots,ds_{2n+1}\in\Omega^1_{K(X)/\CC}$ form a
basis.  Since $2m< 2n$, at least one of the original $ds_1,\ldots, ds_{2n}$ does not lie
in the image of the injection $\Omega^1_{K(Y_1)/\CC}\otimes K(X)\ra\Omega^1_{K(X)/\CC}$.
Say this is $ds_{2n}$.

Now consider the nonzero element $\sigma_1\in R^{2m}(Y_1)$. Up to some factor from $\CC^\times$,
we have $r^*(\sigma_1)=\sigma^m\otimes 1$ inside $R^{2m}(X)=R^{2m}(M)\otimes R^0(E)$, 
where $\sigma\in H^0(M,\Omega^2_M)$ is a symplectic form
on $M$. Since $\sigma$ is nondegenerate, the contraction $\langle\sigma^m, \partial/\partial s_{2m}\rangle\in\Omega^{2m-1}_{K(X)/\CC}$
with the derivation $\partial/\partial s_{2m}$ is nonzero. On the other hand, we have
$\langle r^*(\sigma_1),\partial/\partial s_{2m}\rangle =0$, contradiction.
\qed

\medskip
Recall that a \emph{bielliptic} surface is a minimal surface $S$ with $\kappa=0$ and
$b_2=2$. Equivalently, $S$ has an \'etale covering that is an abelian surface, but is
not an abelian surface itself. Note that the term \emph{hyperelliptic surface} is also frequently used.
We shall discuss such surfaces more thoroughly in Section \ref{bielliptic surfaces}. 

\begin{theorem}
Let $S$ be a bielliptic surface and $Y=\Hilb^n(S)$ for some $n\geq 2$.
Then there is an \'etale covering $\widetilde{Y}\ra Y$ so that $\widetilde{Y}$ is   the product of an elliptic curve
and a  Calabi--Yau manifold of dimension $2n-1$.
\end{theorem}

\proof
Since $h^{1,0}(S)=1$, the Albanese variety of $S$
is an elliptic curve $E$. Let $a:S\ra E$ be the Albanese map. Then $a$ is surjective, and
$\O_E=a_*(\O_S)$, by the universal property. Consider the composite
map
\begin{equation}
\label{hilbert chow}
f:Y=\Hilb^n(S)\lra\Sym^n(S)\lra\Sym^n(E)\lra E,
\end{equation}
where the first arrow is the Hilbert--Chow morphism, the second arrow is
induced from the Albanese map, and the last arrow is given by addition.
Clearly, $f:Y\ra E$ is surjective with $\O_E=f_*(\O_Y)$.
As discussed in Section \ref{bogomolov decomposition}, we have $\pi_1(Y)=H_1(S,\ZZ)$, thus $b_1(Y)=2$ and $h^{1,0}(Y)=1$.
It follows that $f:Y\ra E$ is the Albanese map.

Let $Y_0=f^{-1}(0)$ be the fiber over the origin.
The idea now is to apply Proposition \ref{odd cy} to $Y_0$.
As discussed in Section \ref{bogomolov decomposition} $Y$ has $c_1=0$, thus the same holds for $Y_0$.
Next, we have to verify that $\pi_1(Y_0)$ is finite.
Clearly, the Albanese map is a Serre fibration, such that we have 
an exact sequence
$$
\pi_2(E)\lra\pi_1(Y_0)\lra\pi_1(Y)\lra\pi_1(E)\lra 0.
$$ 
The term on the left vanishes, and the fundamental groups of $Y$ and $E$ are
finitely generated abelian groups of rank two. It follows that $\pi_1(Y_0)$
is finite.

It remains to construct a birational map $X_0\dashrightarrow Y_0$ as in Proposition \ref{odd cy}.
To do so, choose a finite abelian covering $g:A\ra S$ for some abelian surface $A$
and let $X=\Hilb^n(A)$. Then $X$ has a finite covering that is the product of
an abelian surface and the hyperk\"ahler manifold $\Km^n(A)$.
Consider the composite map $A\ra S\ra E$. The kernel of this map 
is a subgroup scheme, whose connected component of the origin $A_0\subset A$ is an
elliptic curve.
Consider the composite map
$$
X=\Hilb^n(A)\lra\Sym^n(A)\lra A.
$$
The fiber $X_0\subset X$ over the
elliptic curve $A_0\subset A$ has $c_1=0$, and admits a finite \'etale covering that
is a product of an elliptic curve and $\Km^n(A)$.
The finite surjection $\Sym^n(A)\ra\Sym^n(S)$ defines a dominant rational map
$\Hilb^n(A)\dashrightarrow\Hilb^n(S)$, and it is easy to see that the latter restricts to 
a dominant rational map $X_0\ra Y_0$.
Thus we may apply Proposition \ref{odd cy} to finish the proof.
\qed

\section{First examples of Enriques manifolds}
\mylabel{first examples}

In this section we shall construct the first examples of higher-dimensional
Enriques manifolds. Suppose   that $S'$ is an Enriques surface, and let $S\ra S'$ be its universal covering.
Then $S$ is a K3 surface endowed with a free action of $G=\ZZ/2\ZZ$
corresponding to a fixed point free involution $\iota:S\ra S$.
This induces a $G$-action on $\Hilb^n(S)$.
In light of Proposition \ref{index divides}, such an action cannot be free for $n$ even.

\begin{proposition}
Suppose $n\geq 1$ is odd. Then the induced $G$-action on $X=\Hilb^n(S)$ is free,
such that $Y=X/G$ is an Enriques manifold of dimension $\dim(Y)=2n$
with index $d=2$.
\end{proposition}

\proof
There is no $\iota$-invariant zero-cycle $\sum_{i=1}^n x_i$ of odd length
on $S$, because the involution $\iota:S\ra S$ is fixed point free.
Thus the induced $G$-action on the symmetric product $\Sym^n(S)$ is free. Since the Hilbert-Chow morphism
$\Hilb^n(S)\ra \Sym^n(S)$ is equivariant, the $G$-action on $\Hilb^n(S)$ must be free as well.
\qed

\medskip
We now turn to Beauville's generalized Kummer varieties.
Fix two elliptic curves $E$ and $E'$,
and consider the abelian surface $A=E\times E'$.
Choose  a point  $a'\in E'$ of order two and an arbitrary point $a\in E$,
and consider the involution
$$
\iota:A\lra A,\quad (b,b')\longmapsto (-b+a,b'+a').
$$
Such maps were studied in connections with cohomologically trivial automorphisms
and are attributed to Lieberman, compare \cite{Mukai; Namikawa 1984}.
The induced action of $G=\ZZ/2\ZZ$ is free, because it is free on 
the second factor. 
Now consider the induced action  on the Hilbert scheme $\Hilb^{n+1}(A)$ and the symmetric product $\Sym^{n+1}(A)$. 
Note  that the addition map $s:\Sym^{n+1}(A)\ra A$ is not equivariant.

\begin{proposition}
Suppose $n\geq 1$ is odd and that $a\in E$ satisfies $(n+1)a=0$ and $(n+1)/2\cdot a\neq 0$.
Then the subset $\Km^n(A)\subset\Hilb^{n+1}(A)$ is $G$-invariant, and the induced
$G$-action on $X=\Km^n(A)$ is free. Whence $Y=X/G$ is an Enriques
manifold of dimension $\dim(Y)=2n$ and index $d=2$.
\end{proposition}

\proof
Let $Z\subset \Sym^n(A)$ be the subscheme of zero-cycles of length $n+1$ on $A$
whose sum in $A$ is the origin $0\in A$. Then we have a commutative diagram
$$
\begin{CD}
X @>>> Z @>>> 0\\
@VVV @VVV @VVV\\
\Hilb^{n+1}(A) @>>> \Sym^{n+1}(A) @>>> A,
\end{CD}
$$
whose squares are cartesian. The Hilbert--Chow morphism is equivariant,
thus it suffices to check that the subset $Z\subset \Sym^{n+1}(A)$ is $G$-invariant
and disjoint from the fixed locus.

Let $\sum_{i=1}^{n+1} x_i$ be a zero-cycle of length $n+1$ on $A$,
and write $x_i=(b_i,b_i')$ with respect to the decomposition $A=E\times E'$.
Suppose that $\sum_{i=1}^{n+1}x_i=0$, where now summation is the actual sum in $A$.
Then $\sum_i b_i=0$ in $E$ and $\sum_ib_i'=0$ in $E'$.
Applying the involution $\iota$ yields a zero-cycle on $A$ summing up to
$$
\sum_i\iota(x_i)=\sum_i(-b_i+a,b_i'+a')=\sum_i(-b_i,b_i') + (n+1)(a,a') = (0,0).
$$
It follows that $Z\subset\Sym^{n+1}(A)$ is $G$-invariant.

Now let $p\in\Sym^{n+1}(A)$ be a $G$-fixed zero-cycle. We check that $p\not\in Z$:
Since $G$ acts freely on $A$, the zero-cycle has the form
$p=\sum_{i=1}^{m} (x_i +\iota(x_i))$
for some closed points $x_1,\ldots,x_m\in A$ where $m=(n+1)/2$.
As above, write $x_i=(b_i,b_i')$ with respect to the decomposition $A=E\times E'$.
Computing the sum   $\sum_{i=1}^{m} (x_i +\iota(x_i))$ in $A$ and projecting onto $E$, we obtain
$\sum_{i=1}^m (b_i-b_i+a) = ma\neq 0$. It follows $p\not\in Z$.
\qed

\medskip
We now introduce the following slightly vague but useful shorthand notation $Y=Q_d X$ in order to refer
to a construction   of Enriques manifolds $Y$  of  index $d\geq 2$   as a quotient of 
a class of hyperk\"ahler manifolds $X$ by some free action of $G=\ZZ/d\ZZ$. For example, we write $Q_2\Hilb^n(S)$ 
for the quotients of Hilbert schemes for K3 surfaces,
and $Q_2\Km^n(E\times E')$ to denote the quotients of Beauville's generalized Kummer
varieties attached to the product of elliptic curve.
We shall generalize these constructions in the next two sections.

\section{Stable sheaves on K3 surfaces}
\mylabel{stable sheaves}

We now generalize the construction $Q_2\Hilb^n(S)$ from the preceding section
using  moduli spaces of stable sheaves.
Throughout this section, $S'$ is an Enriques surface, with universal covering $S\ra S'$, such
that  $S$ is a $K3$-surface endowed with a free action of  $G=\pi_1(S')$,
corresponding to a free involution $\iota:S\ra S$. 

Recall that if $\shF$ is a coherent sheaf on $S$
of rank $r=\rank(\shF)$, first Chern class $l=c_1(\shF)$, and Euler characteristic $\chi=\chi(\shF)$,
then its Mukai vector is
$$
v(\shF)=\ch(\shF)\sqrt{\Todd(S)}=(r,l,\chi-r)\in H^\ev(S,\ZZ).
$$
Let $v\in H^\ev(S,\ZZ)$ be a Mukai vector with $v^2\geq 4$, and $H\in\NS(S)_\RR$ be a polarization.
If   Mukai vector and    polarization are $G$-invariant, then the $G$-action on $S$ induces a $G$-action on
the moduli space $M_H(v)$ of $H$-stable sheaves  on $S$ with Mukai vector $v(\shF)=v$, 
which is a smooth scheme of dimension $v^2+2$ with a symplectic structure but not
necessarily proper.
On the other hand, if $v$ is primitive and $H$ is very general, then $M_H(v)$ is proper, and indeed
a hyperk\"ahler manifold. We have to ensure that the preceding conditions hold simultaneously.

\begin{proposition}
\mylabel{enriques picard}
The following three conditions are equivalent:
\begin{enumerate}
\item The K3 surface $S$ has Picard number $\rho(S)=10$.
\item The canonical map $\Pic(S')\ra\Pic(S)$ is surjective.
\item The $G$-action on $\Pic(S)$ is trivial.
\end{enumerate}
\end{proposition}

\proof
The spectral sequence $H^p(G,H^q(S,\O_S^\times))\Rightarrow H^{p+q}(S',\O_{S'}^\times)$
yields an exact sequence
$$
0\lra H^1(G,\CC^\times)\lra\Pic(S')\lra\Pic(S)^G\lra H^2(G,\CC^\times).
$$
The term on the right vanishes because the group $\CC^\times $ is divisible and $G$ acts trivially.
Furthermore, Enriques surfaces have Picard number $\rho=10$,  and the statement easily follows.
\qed

\begin{proposition}
\mylabel{action moduli}
Suppose $S'$ is an Enriques surface whose corresponding K3 surface $S$ has Picard number
$\rho(S)=10$. Let $v=(r,l,\chi-r)\in H^\ev(S,\ZZ)$ be a Mukai vector and $H\in\NS(S)_\RR$ be
a polarization. Then $v, H$ are $G$-invariant, such that the $G$-action on $S$
induces a $G$-action on $M_H(v)$. If $\chi\in\ZZ$ is odd and $M_H(v)\neq\emptyset$,
then this  $G$-action on $M_H(v)$ is free.
\end{proposition}

\proof
According to Proposition \ref{enriques picard}, the $G$-action on $\Pic(S)$ is trivial.
Consequently  the classes $l,H\in\NS(S)_\RR$ are invariant, hence also the Mukai vector $v\in H^\ev(S,\ZZ)$.
It follows that $\shF\mapsto \iota^*(\shF)$ defines a $G$-action on the moduli space
$M_H(v)$ of $H$-stable sheaves $\shF$ with Mukai vector $v(\shF)=v$.

Now suppose $\chi\in\ZZ$ is odd. Seeking a contradiction, we assume that the $G$-action 
has a fixed point $x\in M_H(v)$. The corresponding coherent sheaf $\shF$  
is then isomorphic to $\iota_*(\shF)$. Choose an isomorphism $h:\shF\ra\iota_*(\shF)$.
Using $\iota_*(\iota_*(\shF))=\shF$, we may regard $\iota_*(h)\circ h$ as an endomorphism of $\shF$.
Being stable, the sheaf $\shF$ is simple, whence $\iota_*(h)\circ h=\lambda\id$ for some
scalar $\lambda\in\CC^\times$. Multiplying $h$ with a square root of $1/\lambda$, we may
assume $\lambda=1$, and this means that $h:\shF\ra\iota_*(\shF)$ defines
a $G$-linearization.

Let $p:S\ra S'$ be the canonical projection. By descend we have $\shF=p^*(\shF')$ for some coherent sheaf $\shF'$ on $S'$. Consequently
$\chi=\chi(\shF)=|G|\chi(\shF')$ is even, contradiction.
\qed

\begin{theorem}
\mylabel{enriques stable}
Suppose $S'$ is an Enriques surface whose corresponding K3 surface has Picard number $\rho(S)=10$.
Let $v=(r,l,\chi-r)\in H^\ev(S,\ZZ)$ be a primitive Mukai vector with $v^2\geq 0$
and $\chi\in\ZZ$ odd. Then for very general polarizations $H\in\NS(S)_\RR$, 
the moduli space $X=M_H(v)$ is a hyperk\"ahler manifold endowed with a free $G$-action, and 
$Y=X/G$ is an Enriques manifold of dimension $v^2+2$ and index $d=2$.
\end{theorem}

\proof
The moduli space $M_H(v)$ is a hyperk\"ahler manifold because $H$ is very general and $\chi$ is primitive
(\cite{O'Grady 1997}, together with \cite{Yoshioka 2001}, Proposition 4.12. Compare also
\cite{Huybrechts; Lehn 1997}, Chapter 6.2).
The $G$-action is free by Proposition \ref{action moduli}.
\qed

\begin{remark}
Under the assumption of the preceding theorem, we have 
$$
2n=\dim M_H(v)=v^2+2=l^2-2r(\chi-r)+2
$$
with $n$ odd,in accordance with Proposition \ref{index divides}.
Indeed,   $l^2\cong 0$ modulo $4$, because any invertible sheaf on $S$ comes from  $S'$,
and $S\ra S'$ has degree two, and the intersection form on $\Pic(S')$ is even. Moreover  
$-2r(\chi-r)\cong 0$ modulo $4$ regardless of $r$, because $\chi$ is odd.
\end{remark}

\medskip
We now check that the conditions of Proposition \ref{enriques picard} indeed hold for very general Enriques surfaces.
For this we need the Global Torelli Theorem for Enriques surfaces, which is due to Horikawa 
(\cite{Horikawa 1978} and \cite{Horikawa 1978}; see also \cite{Barth; Peters; Van de Ven 1984})
Let $L=E_8(-1)^{\oplus 2} \oplus H^{\oplus 3}$
be the K3 lattice, endowed with the   involution 
$$
\iota:L\lra L,\quad (x,y,z_1,z_2,z_3)\longmapsto(y,x,-z_1,z_3,z_2).
$$
Then   each Enriques surface $S'$ admits a marking, that is, an equivariant isometry $L\ra H^2(S,\ZZ)$.
Obviously, the  antiinvariant sublattice $L'\subset L$ is isomorphic to $E_8(-2)\oplus H\oplus H(2)$.
Inside the period domain for K3 surfaces 
$$
D=\left\{[\sigma]\in\PP(L\otimes\CC)\mid (\sigma\cdot\sigma)=0, (\sigma\cdot\bar{\sigma})>0\right\},
$$
the period domain for Enriques surfaces is defined as the intersection
$$
D'=D\cap\PP(L'\otimes\CC).
$$
The discrete group $\Gamma'=\left\{g'\in\Aut(L')\mid \text{  $\exists g\in\Aut(L)$ with $g\iota=\iota g$ and $g'=g|L'$ }\right\}$
acts on $D'$ properly discontinuously, such that the quotient $D'/\Gamma'$ is a normal complex space.
In fact, it acquires the structure of a quasiprojective scheme.
To obtain the coarse moduli space of Enriques surfaces, one has to remove certain divisors.
For each $d\in L'$, consider the divisor $H'_d=\left\{[\sigma]\in D'\mid (\sigma\cdot d)=0\right\}\subset D'$,
and let $H'=\bigcup H_d$ be the union over all $d\in L'$ with $(d\cdot d)=-2$.
It is   known that 
$H'/\Gamma'\subset D'/\Gamma'$ is an irreducible divisor \cite{Allcock 2000}.
The Global Torelli Theorem  asserts that
$$
D'/\Gamma'\smallsetminus H'/\Gamma'=(D'\smallsetminus H')/\Gamma' 
$$ 
is
a coarse moduli space for Enriques surfaces, such that its closed points    bijectively correspond
to isomorphism classes of Enriques surfaces. We note in passing  that the coarse moduli space is quasiaffine \cite{Borcherds 1996}.

\begin{proposition}
The set of    points in the coarse moduli space $(D'\smallsetminus H')/\Gamma'$ corresponding 
to   isomorphism classes of Enriques surfaces $S'$ 
whose  universal covering $S$ has Picard number $\rho(S)=10$ is the complement
of the union of countable many prime divisors.
\end{proposition}

\proof
This condition on marked Enriques surface $S'$ means that no $d\in L'$ is contained
in $\NS(S)\subset H^2(S,\ZZ)$. In other words, each $d\in L'$ must be orthogonal
to $\sigma\in H^{2,0}(S)\subset H^2(S,\CC)$. Consequently, the period $[\sigma]\in D'$
lies in the complement of the union $\bigcup _{d\in L'} H_d'$.
The latter is union is locally finite, because no point in $D'$ is contained
in an intersection $H'_{d_1}\cap\ldots\cap H'_{d_s}$ where $d_1,\ldots,d_s\in L'$ generate
a subgroup of finite index. Therefore, the union $\bigcup_{d\in L'}H'_d\subset D$
is a divisor, which is clearly $\Gamma'$-invariant.
Hence the subset $(\bigcup_{d\in L'}H'_d)/\Gamma'\subset D'/\Gamma'$ is locally
a Weil divisor,  whence itself the countable union of prime divisors.
\qed

\section{Bielliptic surfaces}
\mylabel{bielliptic surfaces}

We now generalize the construction method $Q_2\Km^n(A)$   using the theory of bielliptic
surfaces. Recall that a minimal surfaces $S$ of Kodaira dimension $\kappa=0$ and second Betti number $b_2=2$ is called
\emph{bielliptic}, or \emph{hyperelliptic}. An equivalent condition is that $S$ is not isomorphic to an abelian surface but it admits  a finite \'etale  covering by an abelian surface.
It turns out that the canonical class $\omega_S\in\Pic(S)$ has finite order
$d\in\left\{2,3,4,6\right\}$, and that the corresponding finite \'etale     covering $A\ra S$ is indeed an abelian surface. 
Note that this is an abelian Galois covering whose Galois group $G$ is cyclic of order $d$,
and that 
$$
\pr_*(\O_A)=\O_S\oplus\omega_S\oplus\ldots\oplus\omega_S^{\otimes d-1}.
$$
We call $A\ra S$ the \emph{canonical covering} of $S$.
It turns out that $A$ is isogeneous to a product of elliptic curves.
More precisely, 
there is   a finite \'etale Galois covering $\tilde{A}\ra S$ factoring over $A$, where $\tilde{A}=E\times F$ 
is a product of elliptic curves with  Galois group  of the
form $G\times\tilde{T}$, where $\tilde{T}=\ker(\tilde{A}\ra A)$ is a finite subgroup  and the Galois action
of $G$ and $\tilde{T}$ 
on $\tilde{A}=E\times F$ split into direct product actions.

There are  only seven possibilities, and the whole situation was classified by Bagnera and de Franchis.
We now recall this result in a form we shall use: Consider the complex roots of unity 
$$
i=e^{2\pi i/4},\quad\omega=e^{2\pi i/3},\quad\zeta=e^{2\pi i/6},
$$
 and write $E=\CC/(\ZZ+\tau_1\ZZ)$ and $F=\CC/(\ZZ+\tau_2\ZZ)$, where $\tau_1,\tau_2\in\HH$ are periods.
Also, let $z\in F$ be an  arbitrary element.
The classification is as follows  (compare \cite{Bombieri; Mumford 1977} and \cite{Bennett; Miranda 1990}):

\hspace{-4em}
\begin{array}[t]{|cclll|}
\hline
\tau_2 & d & G\times\tilde{T} & 
\multicolumn{2}{l|}{\text{action of generators on $\tilde{A}=E\times F$}}\\
\hline
\text{arbitrary} & 2 & \ZZ/2\ZZ               & (e,f)\mapsto (e+1/d,-f+z)&\\
\zeta            & 3 & \ZZ/3\ZZ               & (e,f)\mapsto (e+1/d,\omega f+z)&\\
i                & 4 & \ZZ/4\ZZ               & (e,f)\mapsto (e+1/d,if+z)&\\
\zeta            & 6 & \ZZ/6\ZZ               & (e,f)\mapsto (e+1/d,\zeta f+z)&\\
\hline
\text{arbitrary} & 2 & \ZZ/2\ZZ\times\ZZ/2\ZZ & (e,f)\mapsto (e+1/d,-f+z),       & (e,f)\mapsto(e+\tau_1/2,f+1/2)\\
\zeta            & 3 & \ZZ/3\ZZ\times\ZZ/3\ZZ & (e,f)\mapsto (e+1/d,\omega f+z),& (e,f)\mapsto(e+\tau_1/3,f+(1+\zeta)/3)\\
i                & 4 & \ZZ/4\ZZ\times\ZZ/2\ZZ & (e,f)\mapsto (e+1/d,if+z),       & (e,f)\mapsto(e+\tau_1/2, f+(1+i)/2)\\
\hline
\end{array}

\bigskip
So far, the element $z\in F$ is irrelevant, 
because each action is conjugate via a suitable translation
to the ``standard action'' with $z=0$. 
Now fix an integer $n\geq 1$, and let $S$ be a hyperelliptic surface.
The $G$-action on the canonical covering $A$ induces an action on $\Hilb^{n+1}(A)$.
Note that the addition map $\Hilb^{n+1}(A)\ra A$ is not equivariant. Nevertheless, we seek conditions
under which the zero fiber   is invariant, and here our $z\in F$ comes into play:

\begin{proposition}
\mylabel{invariant subset}
Suppose that $d\mid n+1$ and that $z\in F[n+1]$. Then the subset $\Km^n(A)\subset\Hilb^{n+1}(A)$
is $G$-invariant.
\end{proposition}

\proof
Let $\Sym^n_0(A)\subset \Sym^{n+1}(A)$ be the subscheme of zero-cycles
summing up to the origin $0\in A$. Then we have a commutative diagram
$$
\begin{CD}
\Km^n(A) @>>> \Sym^{n+1}_0(A) @>>> 0\\
@VVV @VVV @VVV\\
\Hilb^{n+1}(A) @>>> \Sym^{n+1}(A) @>>> A,
\end{CD}
$$
whose squares are cartesian. The Hilbert--Chow morphism is equivariant,
thus it suffices to check that the subset $\Sym^{n+1}_0(A)\subset\Sym^{n+1}(A)$ is $G$-invariant.

Let $\sum_{i=1}^{n+1} x_i$ be a zero-cycle of length $n+1$ on $A$ summing up to zero.
Choose lifts $\tilde{x}_i\in\tilde{A}$ 
and write $\tilde{x}_i=(e_i,f_i)$ with respect to the decomposition $\tilde{A}=E\times F$.
Let $g\in G$ be the  canonical generator  and write
$g\cdot (e,f)=(e+1/d,\xi f +z)$ with $\xi\in\left\{-1,\omega,i,\zeta\right\}$, as in the Table.
Application of the this automorphism   yields a zero-cycle on $\tilde{A}$ summing up to
$$
\sum_{i=1}^{n+1}g\tilde{x}_i=\sum_i(e_i+1/d,\xi f_i+z)=(1\times \xi)(\sum_i \tilde{x}_i) + (n+1)(1/d,z).
$$
The second summand vanishes by our assumptions, and   $\sum_i\tilde{x}_i$
lies in $\tilde{T}$. It remains to check that $\tilde{T}\subset\tilde{A}$ is invariant under
the automorphism given by $(1\times\xi)$, which is an easy direct computation.
\qed

\medskip
Next, we study  fixed points on Hilbert schemes and symmetric products:

\begin{proposition}
If there is a $G$-fixed point $p\in\Sym^{n+1}(A)$, then $d\mid n+1$.
\end{proposition}

\proof
By induction on $n\geq 1$. Write $p=\sum_{i=1}^{n+1}x_i$.
Since $G$ acts freely on $A$, the $G$-orbit $G\cdot x_1\subset A$ consists
of $d$ pairwise different points, and is contained in the support of $p$, such that $p-Gx_1$ is a $G$-fixed zero cycle
on $A$ of length $n+1-d$. The latter is divisible by $d$ by induction, whence the same holds for $n+1$.
\qed

\medskip
We now make an auxiliary computation: Let $S$ be a bielliptic surface
whose canonical class $\omega_S\in\Pic(S)$ has order $d$,
and consider the  action of $G=\ZZ/d\ZZ$ on  $\tilde{A}=E\times F$. 
Let $g\in G$ be the canonical generator, and write the action as $g(e,f)=(e+1/d,\xi f+z)$
as in the table, with $\xi=-1,\omega,i,\zeta$ for $d=2,3,4,6$, respectively.
Suppose that there is  a $G$-fixed point $p\in\Sym^{n+1}(A)$.
As in the proof for the preceding Proposition, we have
$p=\sum_{i=1}^{m}\sum_{j=0}^{d-1} g^j(x_i)$,
where $m=(n+1)/d$ and $x_1,\ldots,x_m\in A$ are suitable closed points.
Choose lifts $\tilde{x}_i\in\tilde{A}$.

\begin{lemma}
\mylabel{F component}
Assumptions as in the preceding paragraph. Then the
$F$-component of the sum 
$
\sum_{i=1}^m\sum_{j=0}^{d-1}g^j(\tilde{x}_i)\in\tilde{A}=E\times F
$
equals $\sum_{k=1}^{d-1} (d-k)\xi^{k-1}mz\in F$.
\end{lemma}

\proof
Write $\tilde{x_i}=(e_i,f_i)$ with respect to the decomposition $\tilde{A}=E\times F$.
Computing the sum  $\sum_{i=1}^{m}\sum_{j=0}^{d-1} g^j(x_i)$ in $\tilde{A}$ and projecting onto $F$, we obtain
\begin{eqnarray*}
&\sum_{i=1}^m \sum_{j=0}^{d-1}(\xi^jf_i+(\xi^{j-1}+\xi^{j-2}+\ldots +\xi^0)z)\\
=& \sum_{i=1}^m(1+\xi+\ldots+\xi^{d-1})f_i +\sum_{k=1}^{d-1} (d-k)\xi^{k-1}mz.
\end{eqnarray*}
Obviously, the $d$-th root of unity $\xi\neq 1$ is a root of the polynomial $1+T+\ldots+T^{d-1}$,
and the result follows.
\qed

\medskip
We come to the main result of this section. Recall that $\widetilde{T}\subset\widetilde{A}$
is the kernel of $\widetilde{A}\ra A$. Its image under the projection $\widetilde{A}=E\times F\ra F$
is called $T\subset F$.

\begin{theorem}
\mylabel{hyperelliptic enriques}
Suppose $S$ is a hyperelliptic surface whose canonical class has order $d$,
and $A\ra S$ be its canonical covering.
Let $n\geq 1$ be an integer with $d\mid n+1$, and $z\in F[n+1]$.
Write $m=(n+1)/d$ and assume:
\begin{enumerate}
\item If $d=2$  then $mz\not\in  T$.
\item If $d=3$  then $T=0$ and $mz\not\in\ZZ(1+\zeta)/3$.
\item If $d=4$  then $T=0$ and $2mz\not\in\ZZ(1+i)/2$.
\end{enumerate}
Then the subset $\Km^n(A)\subset\Hilb^{n+1}(A)$ is $G$-invariant,
and the induced $G$-action on $X=\Km^n(A)$ is free, such
that $Y=X/G$ is an Enriques manifold of dimension $\dim(Y)=2n$ and index $d$.
\end{theorem}

\proof
We already saw in Proposition \ref{invariant subset} that $\Km^n(A)\subset\Hilb^{n+1}(A)$ is
invariant. Seeking a contradiction, we suppose that the induced action on $X=\Km^n(A)$ is not free.

Let us first consider the cases $d=2$ and $d=3$. Then there is a fixed point on $X$,
and its image under the Hilbert--Chow morphism is a fixed point $p\in\Sym^{n+1}_0(A)$.
Since $G$ acts freely on $A$, we may write $p=\sum_{i=1}^m\sum_{j=0}^{d-1}g^j(x_i)$,
where $x_1,\ldots,x_m\in A$ are closed points, and $g\in G$ is the canonical generator.
Choose lifts $\tilde{x}_i\in\tilde{A}$, and write $\tilde{x_i}=(e_i,f_i)$
with respect to the decomposition $\tilde{A}=E\times F$.
We now use  Lemma \ref{F component} to compute  the $F$-component of the sum $\sum_{i=1}^m\sum_{j=0}^{d-1}g^j(\tilde{x}_i)\in\tilde{A}$:

In case $d=2$, the $F$-component is given by $mz\in F$, which is not contained
in $T\subset F$ by assumption. Whence $\sum_{i=1}^m\sum_{j=0}^{d-1}g^j(\tilde{x}_i)$
is not contained in the kernel $\tilde{T}$ for the homomorphism $\tilde{A}\ra A$.
On the other hand, we have $\sum_{i=1}^m\sum_{j=0}^{d-1}g^j(x_i)=0$ in $A$, contradiction.

Now suppose $d=3$. Then the $F$-component is given by $(2+\omega)mz$.
One easily computes that $(1+\zeta)/3$ generates the kernel of $(2+\omega)$
viewed as an endomorphism of $F[3]$. So our assumption ensures
that $(2+\omega)mz\neq 0$, and we obtain a contradiction as above.

It remains to treat the cases $d=4$. Choose  a point $p\in\Sym^{n+1}_0(A)$ whose
stabilizer is nonzero. 
In case $d=4$, this point is fixed by the unique subgroup $G'\subset G$
of order two. Applying the preceding paragraph to the hyperelliptic surface
$S'=A/G'$, we obtain a contradiction.
\qed

\medskip
In all cases, the element $mz\in F[d]$ or   suitable multiples have
to avoid a 1-dimensional vector subspace in a 2-dimensional
vector space over  certain finite fields.
This can always be done, so the cases are indeed nonvacuous.
Thus:

\begin{theorem}
There are Enriques manifolds of index $d=2,3,4$.
\end{theorem}

\begin{remark}
As pointed out by Sarti and Boissi\'ere, the case $d=6$ seems impossible here.
\end{remark}



\begin{thebibliography}{ccccc}

\bibitem{Allcock 2000}
D.\ Allcock:
The period lattice for Enriques surfaces.
Math.\ Ann.\  317  (2000),  483--488.

\bibitem{Barth; Peters; Van de Ven 1984}
W.~Barth, C.~Peters, A.~Van de Ven:
Compact complex surfaces.
Ergeb.\ Math.\  Grenzgebiete (3) 4,
Springer, Berlin, 1984.

\bibitem{Beauville 1983}
A.\ Beauville:
Vari\'et\'es K\"ahleriennes dont la premi\`ere classe de Chern est nulle.
J.\ Differential Geom.\ 18 (1983) 755--782.

\bibitem{Beauville; Donagi 1985}
A.\ Beauville, R\ Donagi:
La vari\'et\'e des droites d'une hypersurface cubique de dimension $4$. 
C.\ R.\ Acad.\ Sci.\ Paris S\'er.\ I Math.\  301  (1985),   703--706. 

\bibitem{Bennett; Miranda 1990}
C.\ Bennett, R.\  Miranda:
The automorphism groups of the hyperelliptic surfaces.
Rocky Mountain J.\ Math.\ 20 (1990), 31--37. 

\bibitem{Bogomolov 1974a}
F.\ Bogomolov:
K\"ahler manifolds with trivial canonical class.
Izv.\ Akad.\ Nauk SSSR Ser.\ Mat.\ 38 (1974), 11--21.

\bibitem{Bogomolov 1974b}
F.\ Bogomolov:
The decomposition of K\"ahler manifolds with a trivial canonical class. 
Mat.\ Sb.\ (N.S.) 93 (135) (1974), 573--575, 630.

\bibitem{Bombieri; Mumford 1977}
E.~Bombieri, D.~Mumford:
Enriques' classification of surfaces in char.\ p.  II.
In: W.~Baily, T.~Shioda (eds.),
Complex analysis and algebraic geometry, pp.\ 23--42.
Cambridge University Press, London, 1977.

\bibitem{Borcherds 1996}
R.\ Borcherds:
The moduli space of Enriques surfaces and the fake Monster Lie superalgebra.
Topology 35 (1996), 699--710. 

\bibitem{Brion; Kumar 2005}
M.\ Brion, S.\ Kumar:
Frobenius splitting methods in geometry and representation theory.
Progress in Mathematics 231.
Birkh\"auser, Boston, 2005.

\bibitem{Debarre; Voisin 2009}
O.\ Debarre, C.\ Voisin:
Hyper-K\"ahler fourfolds and Grassmann geometry.
Preprint (2009), arXiv:math.AG/0904.3974.

\bibitem{Griffiths; Harris 1978}
P.\ Griffiths, J.\ Harris:
Principles of algebraic geometry.
Wiley-Interscience, New York, 1978.

\bibitem{SGA 1}
A.~Grothendieck et al.:
Rev\^etements \'etales et groupe fondamental.
Lect.\ Notes Math.\  224,
Springer, Berlin, 1971.

\bibitem{Horikawa 1978}
E.\ Horikawa:
On the periods of Enriques surfaces. I.
Math.\ Ann.\ 234 (1978),   73--88.

\bibitem{Horikawa 1978b}
E.\ Horikawa:
On the periods of Enriques surfaces. II.
Math.\ Ann.\ 235 (1978), 217--246. 

\bibitem{Huybrechts; Lehn 1997}
D.\ Huybrechts, M.\ Lehn:
The geometry of moduli spaces of sheaves. 
Aspects of Mathematics E31. 
Vieweg, Braunschweig, 1997.

\bibitem{Huybrechts 1999}
D.\ Huybrechts:
Compact hyper-K\"ahler manifolds: basic results.
Invent.\ Math.\  135  (1999),  63--113.

\bibitem{Iliev; Ranestad 2001}
A.\ Iliev, K.\ Ranestad:
$K3$ surfaces of genus 8 and varieties of sums of powers of cubic fourfolds.  
Trans.\ Amer.\ Math.\ Soc.\  353  (2001),   1455--1468.

\bibitem{Iversen 1970}
B.\ Iversen:
Linear determinants with applications to the Picard scheme of a family of algebraic curves.
Lecture Notes in Mathematics 174.
Springer, Berlin, 1970.

\bibitem{Kollar 1993}
J.\ Koll\'ar:
Shafarevich maps and plurigenera of algebraic varieties.  
Invent.\ Math.\  113  (1993),   177--215.

\bibitem{Lange 2001}
H\ Lange:
Hyperelliptic varieties.  
Tohoku Math.\ J.\  53  (2001),  491--510.

\bibitem{Mukai 1984}
S.\ Mukai:
Symplectic structure of the moduli space of sheaves on an abelian or K3 surface.
Invent.\ Math.\ 77 (1984), 101-116.

\bibitem{Mukai 1987}
S.\ Mukai:
On the moduli space of bundles on K3 surfaces. I.
In: M.\ Atiyah et al. (eds.),
Vector bundles on algebraic varieties.
Tata Inst.\ Fundam.\ Res.\ 11 (1987), 341--413.

\bibitem{Mukai; Namikawa 1984}
S.\ Mukai, Y.\ Namikawa:
Automorphisms of Enriques surfaces which act trivially on the cohomology groups.
Invent.\ Math.\ 77 (1984),  383--397. 

\bibitem{Nikulin 1980}
V.\ Nikulin:
Finite automorphism groups of K\"ahler K3 surfaces. 
Trans.\ Mosc.\ Math.\ Soc.\ 2 (1980), 71--135. 

\bibitem{O'Grady 1997}
K.\ O'Grady:
The weight-two Hodge structure of moduli spaces of sheaves on a $K3$ surface.
J.\ Algebraic Geom.\ 6 (1997),  599--644. 

\bibitem{O'Grady 1999}
K.\ O'Grady:
Desingularized moduli spaces of sheaves on a $K3$.  
J.\ Reine Angew.\ Math.\  512  (1999), 49--117.

\bibitem{O'Grady 2003}
K.\ O'Grady:
A new six-dimensional irreducible symplectic variety.  
J.\ Algebraic Geom.\  12  (2003),  435--505.

\bibitem{O'Grady 2006}
K.\ O'Grady:
Irreducible symplectic 4-folds and Eisenbud-Popescu-Walter sextics. 
Duke Math.\ J.\ 134 (2006), 99--137. 

\bibitem{Yoshioka 2001}
K.\ Yoshioka:
Moduli spaces of stable sheaves on abelian surfaces.  
Math.\ Ann.\  321  (2001),  817--884.
\end{thebibliography}
\end{document}